\documentclass[11pt]{article}
\usepackage{amsfonts,amssymb,amsmath,amsgen}

{\it}{\rm}
{\it}{\rm}
\newtheorem{corol}{COROLLARY}{\bf}{\it}
{\bf}{\it}
{\it}{\rm}
{\bf}{\rm}
\newtheorem{lemma}{LEMMA}[section]{\bf}{\it}
{\it}{\rm}
{\bf}{\rm}
{\it}{\rm}
\newtheorem{prop}{PROPOSITION}{\bf}{\it}
{\it}{\rm}
{\bf}{\rm}
\newtheorem{remark}{Remark}{\bf}{\rm}
\newtheorem{theorem}{THEOREM}{\bf}{\it}
{\bf}{\rm}

\makeatother
\newcommand{\sphlap}{\displaystyle{\not\!\!\Delta}}
\newcommand{\sphgrad}{\displaystyle{\not\!\nabla}}

\newcommand{\ds}{\displaystyle}

\newcommand{\la}{\lambda}

\newcommand{\nab}{\nabla}

\newcommand{\p}{\partial}

\newcommand{\beq}{\begin{equation}}
\newcommand{\eeq}{\end{equation}}
\newcommand{\bna}{\begin{eqnarray}}
\newcommand{\ena}{\end{eqnarray}}
\newcommand{\ba}{\begin{eqnarray*}}
\newcommand{\ea}{\end{eqnarray*}}
\newcommand{\ben}{\begin{enumerate}}
\newcommand{\een}{\end{enumerate}}
\newcommand{\bi}{\begin{itemize}}
\newcommand{\ei}{\end{itemize}}

\newcommand{\Hi}{\dot{\cal H}}

\newcommand{\R}{{\cal R}}
\newcommand{\RR}{{\mathbb R}}

\newcommand{\CC}{{\mathbb C}}

\newcommand{\BS}{{\mathbb S}}

\newcommand{\qed}{\hfill\rule{2mm}{2mm}}

\newcommand{\real}{\mathop{\mathrm{Re}}\nolimits}

\begin{document}

\title{
Strichartz estimates\\ for the Wave and Schr\"odinger  Equations\\
with Potentials of Critical Decay }

\author{
Nicolas Burq
\and
Fabrice Planchon
\and
John G. Stalker
\and
A. Shadi Tahvildar-Zadeh
}

\date{\today}

\maketitle
\begin{abstract}
We prove weighted $L^2$  estimates for the solutions of linear Schr\"odinger and wave
equation with potentials that decay like $|x|^{-2}$ for large $x$, by deducing them
 from estimates on the resolvent of the associated elliptic
operator. We then deduce Strichartz estimates for these equations.
\end{abstract}
\section{Introduction and Main Results}

Consider the following linear equations
\bna\label{eq:lse}
&&i\partial_t u +\Delta u- V(x) u = 0\qquad
u(0,x) = f(x)
\\
\label{eq:lwe}
&&-\partial_t^2 u +\Delta u- V(x) u = 0\qquad
u(0,x) = f(x),\ \p_t u(0,x) = g(x)
\ena
where $\Delta $ is the $n$ dimensional Laplacian. Throughout this paper
we will assume $n\geq 3$.

In \cite{BPST1} we showed that in the case where $V(x)=\frac{a}{|x|^2}$ the solution to the above equations
satisfies generalized spacetime Strichartz  estimates as long
as $a> -(n-2)^2/4$.  We intend to  extend this result to
potentials which, in a sense to be made precise below, behave like the inverse square potential.

Let $\lambda(n)$ be defined as follows
\begin{equation}\label{def:la}
\lambda := \frac{n-2}{2}
\end{equation}
We  also define multiplication operators $\Omega^s$ by
\begin{displaymath}
(\Omega^s \phi)(x) = |x|^s \phi (x).
\end{displaymath}
Let  $\sphlap$ denote
the spherical Laplacian and $\sphgrad$ the spherical gradient
 on the unit sphere. Let $r(x) := |x|$ denote
 the polar radius. For a given function $V\in C^1(\RR^n\setminus\{0\})$
 let $\tilde{V}$ be defined by
\[
\tilde{V}(x) :=  - \partial_r (rV(x)),
\]

We denote the positive and negative parts of a function $V$ by $V_+ :=
\max\{V,0\}$ and $V_- := \max\{-V,0\}$ respectively. Thus $V = V_+ - V_-
$.

In this paper we consider time-independent potentials $V(x) \in C^1(\RR^n\setminus\{0\})$ satisfying
the following assumptions.
\begin{itemize}
\item[{\bf (A1)}]
$\gamma_\pm^2 := \sup_{x\in \RR^n} |x|^2 V_\pm(x) < \infty$
\item[{\bf (A2)}] The operator
$-\sphlap + \Omega^2 V +\lambda^2$ is positive on every sphere, i.e.,
there exists a $\delta>0$ such that for every $r>0$,
\begin{equation}\label{A2}
 \int_{|x| = r} |\sphgrad u(x)|^2  +
 (\lambda^2 + |x|^2V(x) )|u(x)|^2 d\sigma(x)  \geq \delta^2 \int_{|x| = r}
|u(x)|^2 d\sigma(x)
\end{equation}
 \item[{\bf (A3)}] The operator $-\sphlap +\Omega^2 \tilde{V} + \lambda^2$ is
positive on every sphere, i.e. (\ref{A2}) holds with $\tilde{V}$ in
place of $V$.
\end{itemize}
\begin{remark} The potential is thus allowed to have one point
singularity, which without loss of generality we take to be at the
origin of coordinates.  Note that no sign condition is assumed on $V$,
and that only differentiability in the radial direction is actually
used.

Note also that for an inverse-square potential $V=a|x|^{-2}$ assumptions {\bf (A2)} and {\bf (A3)} are the same, and
require that $a>-\lambda^2$.  More generally, for potentials that are homogeneous functions of degree $-2$, i.e.
 $V(x) = |x|^{-2} a(x/|x|)$ with $a$ a function defined on the unit
sphere, these two assumptions are again the same, namely that $-
\sphlap + a +\lambda^2$ be a positive operator on  the unit sphere.
In section 4 we will consider an application where such a potential
appears.
\end{remark}
\begin{remark}
While the approach we present recovers the results of \cite{BPST1} as a
special case, it should be noted that it turns out to be much simpler,
and hence more flexible. In particular, it should be possible to
include time-dependent potentials as well, a strategy which will be
pursued elsewhere.
\end{remark}
Before stating the main results of the paper, let us examine the above
assumptions more closely.  Let $Q(u)$ denote the quadratic form
naturally
associated with the operator $-\Delta+ V$, i.e.
\begin{equation}\label{def:Q}
Q(u) :=  \int_{\RR^n} |\nab u(x)|^2 + V(x)|u(x)|^2 \ dx
\end{equation}
We then have
\begin{prop}\label{prop:equiv}
Under the assumptions {\bf (A1-A2)} for $V$, there are constants $0<c_1\leq 1 \leq c_2$ such that
\[
c_1\| \nab u\|_{L^2}^2 \leq Q(u) \leq c_2 \|\nab u\|_{L^2}^2
\]
\end{prop}

{\em Proof:} Recall the celebrated Hardy's inequality:
\[
\| \Omega^{-1} u \|_{L^2(\RR^n)} \leq \frac{2}{n-2} \| \frac{\partial u}{\partial r} \|_{L^2(\RR^n)}
\]
for $n\geq 3$ (see \cite{KalSchWalWus75} for a proof). By {\bf (A1)}
we thus have
\[
Q(u) \leq \int |\nab u|^2 + \frac{\gamma_+^2}{|x|^2}|u|^2\ dx \leq ( 1 + \frac{\gamma_+^2}{\lambda^2}) \|\nab u \|_{L^2}^2
\]
so that $c_2 =1 + \ds\frac{\gamma_+^2}{\lambda^2}$. On the other hand, by {\bf (A2)},
\begin{eqnarray}
Q(u) & = &  \int_0^\infty \int_{|x| = r} |\frac{\partial u}{\partial |x|}|^2 +
\frac{1}{|x|^2} |\sphgrad u |^2 + V(x) |u(x)|^2 \ d\sigma(x) dr \nonumber\\
& \geq &  \int_0^\infty \frac{1}{r^2} \int_{|x| = r} |\sphgrad u|^2
+ ( r^2 V(x) + \lambda^2) |u(x)|^2 d\sigma(x) dr \nonumber\\
& \geq &  \int_0^\infty \frac{\delta^2}{r^2} \int_{|x| = r}
|u(x)|^2 d\sigma(x) dr = \delta^2 \| \Omega^{-1} u \|_{L^2}^2,\label{PHardy}
\end{eqnarray}
and thus, if we set $c_1 :=\ds \frac{\delta^2}{\delta^2+\gamma_-^2}$, then
using {\bf (A1)} we have
\[
Q(u) - c_1\|\nab u\|_{L^2}^2 \geq \int_{\RR^n} (- c_1 |x|^2 V_-(x) +
(1-c_1)\delta^2)\frac{|u|^2}{|x|^2} dx \geq 0
\]

An important consequence of the above proposition is the following
equivalence result:
\begin{corol}\label{corol:equiv}
 Let $\Hi^s(\RR^n)$ denote the scale of
homogeneous Sobolev spaces based on the powers of the operator $P = -
\Delta + V$.  I.e., the completion of $C^\infty_c(\RR^n\setminus\{0\})$
with respect to the seminorm
\[
\|u \|_{\Hi^s(\RR^n)} := \| P^{s/2} u \|_{L^2(\RR^n)}
\]
If $V$ satisfies {\bf (A1-A2)} then the spaces $\Hi^s$ are
equivalent to the standard Sobolev spaces $\dot{H}^s$ (based on the
powers of $-\Delta$) for $|s| \leq 1$.

\end{corol}

For $s=1$ this follows immediately from the above Proposition, noting that $Q(u) = \|P^{1/2} u\|_{L^2}^2$.
The case $s=-1$  then follows by duality, and by interpolation we get the $s$ in between.

Our main result for the Schr\"odinger equation (\ref{eq:lse}) is
\begin{theorem}\label{thm:strich}
  Let $f \in L^2$, $p\geq 2,q$ such that
  \begin{equation}\label{schrodadmis}
\frac{2}{p}+\frac{n}{q}=\frac{n}{2},\qquad n\geq 3.
\end{equation}
Let $u$ be the unique solution of (\ref{eq:lse}).  Then
there exists a constant $C>0$  such that
  \begin{equation}
    \label{eq:strichartz}
    \|u\|_{L^p_t(L^q_x)} \leq C \|f\|_{L^2}.
  \end{equation}
\end{theorem}
For the wave equation (\ref{eq:lwe}) we have
\begin{theorem}\label{thm:strichwave}
  Let $u$ be the solution to {\rm(\ref{eq:lwe})} with Cauchy data $(f,g)
  \in \dot{H}^{\frac{1}{2}}\times \dot{H}^{-\frac{1}{2}}$.  Let  $p > 2$, and $q$ be
  such that
  $\frac{2}{p}+\frac{n-1}{q}=\frac{n-1}{2}$. Then
  \begin{equation}
    \label{eq:strichartzwave}
    \| (-\Delta)^{\sigma/2} u\|_{L^p_t(L^q_x)} \leq
C (\|f\|_{\dot{H}^{\frac{1}{2}}}+\|g\|_{\dot{H}^{-\frac{1}{2}}}).
  \end{equation}
where $\sigma=\frac{1}{p}+\frac{n}{q}-\frac{n-1}{2}$ (gap condition).
\end{theorem}
\begin{remark}
  Notice we include the end-point for the Schr\"odinger equation while
  we exclude it for the wave equation. This is strictly intended to
  make the argument shorter, and one could adapt the argument from
  \cite{BPST1} to recover the endpoint case for the wave equation as well.
\end{remark}
The strategy for proving the above Strichartz estimates is to deduce
them from the corresponding estimates for the free case $V\equiv 0$,
using Duhamel's principle.  This was the approach taken in
\cite{RodSch} where Strichartz estimates for 3D Schr\"odinger where proved for
(\ref{eq:lse}) under the assumption that $V(x)$ decays like
$|x|^{-2-\epsilon}$ for large $x$. For the 3D wave equation,
dispersive estimates were recently proven in \cite{GeoVis03},
under the following assumptions: $V(x)\geq 0$, $V\in
C^{0^+}(\R\setminus \{0\})$, $V(x)\lesssim \inf
(|x|^{-2+0^+},|x|^{-2-0^+})$ and the usual spectral assumption that
zero is neither an eigenvalue nor a resonance. The method involves
rather delicate resolvent estimates, and once the dispersive estimate
is proven, Strichartz estimates follows by standard considerations. Note that the
$1/|x|^2$ is barely missed and therefore appears like a borderline
case. Indeed, when the potential is admitted to have a singularity at
$x=0$ slightly stronger, space-time estimates may fail, as the example
provided in \cite{Thomas} shows.

Here we follow the strategy from \cite{RodSch} and bypass dispersive
estimates to obtain directly Strichartz estimates. The key ingredient in this
argument is the availability of a weighted spacetime $L^2$ estimate
for the solutions of the above equations.  More precisely, for
(\ref{eq:lse}) one needs
the estimate
\begin{equation}\label{est:KY}
\| \Omega^{-1} u \|_{L^2_tL^2_x} \leq C \|  f \|_{L^2}
\end{equation}
which for the free case is a particular instance of the Kato-Yajima
smoothing estimate \cite{KatYaj89}, while for (\ref{eq:lwe}) the
corresponding estimate needed turns out to be
\begin{equation}\label{est:mor}
\| \Omega^{-1} u \|_{L^2_tL^2_x} \leq C \left( \| f \|_{\dot{H}^{1/2}} + \|
g\|_{\dot{H}^{-1/2}}\right),
\end{equation}
which can be thought of as  a generalization of the Morawetz estimate
\cite{Mor68} (See \cite{Hos97} for a proof of (\ref{est:mor}) in the free case).

Using an abstract machinery largely due to Kato \cite{Kat66} (see \cite[\S XIII.7]{ReeSimIV}),
the above weighted-$L^2$ estimates are deduced from a weighted resolvent
estimate for the elliptic operator $P$ which is a particular self-adjoint extension of
$-\Delta +V$  (see Theorem~\ref{thm:kato}).
We use the method of multipliers to prove this resolvent estimate (see Theorem~\ref{thm:ARA*}).
(see \cite{Per99} where multipliers
are used to prove a similar estimate).

The outline of this paper is as follows:  In Section 2 we prove the
resolvent estimate (\ref{est:res}).  Weighted-$L^2$ and Strichartz
estimates are proved in Section 3.  We consider an application in
Section 4.

\section{Resolvent estimates}

We prove weighted $L^2$ estimates for the resolvent of
$-\Delta + V$.  Note that for the  potentials that we are considering here, $-\Delta+V$ is not  a compact
perturbation of $-\Delta$.
In order to define the resolvent however, we first need to consider self-adjoint
extensions of $-\Delta + V(x)$, which is a symmetric operator but a priori only defined on $C^2(\RR^n\setminus\{0\})$.
  We refer to \cite{PST1} for a similar
discussion in the case of inverse-square potentials $V = a|x|^{-2}$.  In that case
it is well-known that self-adjoint extensions are not unique when $-\lambda^2<a<1-
\lambda^2$. In particular, there are two extensions that are both
rotation and dilation invariant.  One of the two corresponding  domains contains functions
with infinite energy (i.e. infinite $\dot{H}^1$ norm). It will be clear in what follows that
having finite energy is crucial to  the arguments that we present, and that is
why we are going to consider the {\em Friedrichs extension} of the operator $-\Delta +V$, i.e.
using the corresponding quadratic form (\ref{def:Q}) to define the extension. It was
shown in \cite{KalSchWalWus75} (Theorem 3) that for the class of
potentials we are considering, this extension has the property
that the domain of the extended operator is contained in $\dot{H}^1$.  We denote by $P$
the Friedrichs extension of $-\Delta+V$.  $P$ is thus self-adjoint, and
an application of Hardy's inequality, together with
assumption {\bf (A2)} implies that $P$ is a positive operator, and
\[
\sigma(P) = \sigma_{ac}(P) = [0,+\infty).
\]
 It  follows that the resolvent  $R(\mu) := (P-\mu)^{-1}$ is a well-
defined bounded operator on $L^2$ for $\mu \not\in \RR^{+}$.   The goal of this section is to prove
\begin{theorem}\label{thm:ARA*}
Let $V(x)$ satisfy {\bf(A1-A3)}.  Then there is a constant $C>0$ such that
\begin{equation}\label{est:res}
\sup_{\mu \notin \RR^{+}} \| \Omega^{-1} R(\mu)\Omega^{-1} f\|_{L^2} \leq C \|  f\|_{L^2}
\end{equation}
\end{theorem}

The  proof uses the method of multipliers, and is based on
 Morawetz's radial identity~\cite{Mor75}. Given $f \in L^2(\RR^n)$ and $\mu \in \CC\setminus \RR^+$, let $u\in D(P)\subset H^1_0(\RR^n)$ be the
unique solution of the inhomogeneous Helmholtz equation
\begin{equation}\label{eq:helmf}
P u + z^2 u = f
\end{equation}
where $z = \sqrt{-\mu}$, with the branch chosen such that $\real z=\sigma >0$. Thus $u = R(\mu)f$.
In order to carry out the integration by parts argument below, one needs to know something about
the behavior at the origin and at infinity of
$u$, to check that the contributions of these points have ``the good sign''. In the case of the potential $a|x|^{-2}$, this can be done by using the explicit asymptotic behavior of Hankel functions near $0$ and infinity. In the general case, it is actually enough to know that $u \in {H}^1(\RR^n)$ (but the argument requires some care, see below).

To prove (\ref{est:res}) we first note that by density, we can take $f\in C^\infty_0(\RR^n\setminus\{0\})$.
 Then $u$ is clearly a classical solution of (\ref{eq:helmf}).  Let $w:\RR^+\times\BS^{n-1} \to \CC$ be defined
by
\[
w(r,\theta) := r^{\lambda+1/2} e^{rz} u(r\theta).
\]
Then $w$ satisfies
\begin{equation}\label{eq:w}
-\partial_{r}^2 w -\frac{1}{r^2}\sphlap w + (\lambda^2-\frac{1}{4}+r^2V(r\theta))\frac{w}{r^2} + 2z\partial_r w
 =e^{ z r}r^{\la+1/2}f
\end{equation}
For $R>\epsilon>0$ fixed, let $\phi = \phi_{\epsilon,R}(r)$ be a smooth cut-off function, $0\leq \phi\leq 1$,
 that is zero outside $[0,R+1]$ and is equal to
one on $[\epsilon,R]$.
Multiplying (\ref{eq:w}) by $$re^{-2r\sigma}\phi(r)\partial_r\bar{w}$$ and taking the real part, we obtain
\begin{eqnarray*}
\lefteqn{-\frac{1}{2} r e^{-2r\sigma}\phi \partial_r |\partial_r w|^2 +
 \frac{1}{2r}e^{-2r\sigma}\phi \partial_r |\not\!\nabla w|^2 +
2\sigma r e^{-2r\sigma} \phi  |\partial_r w|^2 }\\
&& + \frac{1}{2r}e^{-2r\sigma} \phi (\la^2-\frac{1}{4}+r^2 V)\partial_r|w|^2
 = \real(r^{\la+3/2}\phi e^{r(z-2\sigma)} \partial_r\bar{w} f)
\end{eqnarray*}
Integrating the above equality on $\RR^+\times \BS^{n-1}$ and performing the integration by parts we obtain
\begin{eqnarray}\label{ibp}
\lefteqn{\frac{1}{2} \int_0^\infty\int_{|\theta|=1}
\phi e^{-2r\sigma} \left[ (1+2r\sigma) (|\partial_r w|^2 -\frac{1}{4r^2}|w|^2)\right.}\nonumber\\
&& + \left. \frac{1}{r^2}\left( |\sphgrad w|^2 + (r^2 \tilde{V}(r\theta)+\lambda^2)|w|^2 \right)\right.\nonumber\\
&& + \left. \frac{2r\sigma}{r^2}\left( |\not\!\nabla w|^2 + (r^2 V(r\theta) + \lambda^2)|w|^2 \right) \right] d\theta dr \nonumber \\
& + &  \frac{1}{2} \int_0^\infty\int_{|\theta| = 1} r e^{-2r\sigma}\phi'(r) \left[ |\partial_r w|^2 +\frac{1}{4r^2} |w|^2\right.\nonumber\\
&& - \left. \frac{1}{r^2}
\left( |\sphgrad w |^2 + (r^2 V(r\theta)+\lambda^2) |w|^2 \right) \right] d\theta dr \nonumber\\
& = & \int_0^\infty \int_{|\theta| =1} \real (r^{\la+3/2}\phi e^{r(z-2\sigma)} \partial_r\bar{w} f) d\theta dr
\end{eqnarray}
By Cauchy's inequality, for any $a>0$ the right hand side of the above is less than or equal to
\begin{equation}\label{peterpaul}
\frac{1}{4a^2} \|  \Omega f \|_{L^2}^2 + a^2 \int\int \phi e^{-2r\sigma} |\partial_r w|^2 d\theta dr.
\end{equation}
The difficulty in the analysis of~\eqref{ibp} is twofold: first we have to check that the first integral in the left-hand side controls $w$ in some suitable space. Second we have to show that the contributions of the second integral in this left-hand side are non negative as $\varepsilon\rightarrow 0$ and $R \rightarrow + \infty$. We start by considering the second problem.
We note that $\mbox{supp}\phi'\subset I_\epsilon \cup I_R$, where $I_\epsilon :=[0,\epsilon]$ and $I_R :=[R,R+1]$.
On $I_\epsilon$ we have $0\leq \phi'\leq C/\epsilon$
and on $I_R$ we know that $-C\leq \phi'\leq 0$.  Since the left hand side of (\ref{ibp}) is to be estimated from below, we only
need to estimate the negative terms in this integral. In particular, it is enough to show
\begin{eqnarray}
\lim_{R \to \infty} \int_R^{R+1} \int_{|\theta| = 1} r e^{-2r\sigma}(|\partial_r w|^2 + \frac{1}{4r^2} |w|^2) d\theta dr = 0 \label{infty}\\
\lim_{\epsilon \to 0} \int_0^\epsilon \int_{|\theta| = 1} \frac{1}{r^2}( |\sphgrad w^2| + |w|^2)d\theta dr = 0\label{zero}
\end{eqnarray}

Let us first consider (\ref{infty}).  It is in fact
enough to show that there exists a sequence $R_n \to \infty$ along which it
holds. We note that
\[
\partial_r w = r^{\lambda+1/2} e^{rz} (( \partial_r  + z) u + \frac{\lambda+1/2}{r}u).
\]
We thus have, using Hardy's inequality, that
\begin{eqnarray*}
\lefteqn{ \int_R^{R+1} \int_{|\theta| = 1} re^{-2r\sigma}( |\partial_r w|^2 + \frac{1}{4r^2} |w|^2) d\theta dr
}\hspace{1in}\\
&\leq& C(n) \int_R^{R+1} \int_{|\theta| = 1} r( |\partial_r u|^2 + |z|^2|u|^2) d\theta r^{n-1}dr \\
& \leq & C(n,|z|) \int_R^{R+1} r h(r) dr
\end{eqnarray*}
where $$ h(r) := r^{n-1}\int_{|\theta|=1} |\partial_r u(r\theta)|^2 + |u(r\theta)|^2 d\theta$$
so that by virtue of $u$ being in $H^1(\RR^n)$, we know
 $$ \int_0^\infty h(r) dr <\infty.$$
  It thus follows that given
$\mu_m >0$ there exists a sequence $R^{(m)}_n \to \infty $ such that
 $$ \int_{R^{(m)}_n}^{R^{(m)}_n+1} h(r) dr < \frac{\mu_m}{R^{(m)}_n},$$
because otherwise the integral $\int_0^\infty h dr$ would diverge.  Using a diagonal argument it thus follows that there
exists a sequence $R_n \to \infty$ such that
$$\int_{R_n}^{R_n+1} r h(r) dr \to 0\qquad \mbox{ as } n \to \infty,$$
which establishes (\ref{infty}) along a sequence.

Similarly, using that the $H^1$ norm of $u$ on a ball is finite, we have
\[
\int_0^\epsilon \int_{|\theta| = 1} \frac{1}{r^2}( |\sphgrad w^2| + |w|^2)d\theta dr \leq
 C \int_0^\epsilon \int_{|\theta| = 1} |\nabla u|^2 r^{n-1}d\theta dr \to 0
\]
as $\epsilon \to 0$,
establishing (\ref{zero}). \par
 We can thus focus our attention on the first integral on the left in (\ref{ibp}). Using the assumptions {\bf (A2), (A3)} on the potential, it
can be estimated from below by
\[
\frac{1}{2} \int_0^\infty\int_{|\theta|=1}
\phi e^{-2r\sigma}  (1+2r\sigma) \left[ |\partial_r w|^2 + (\delta^2 - \frac{1}{4})\frac{|w|^2}{r^2} \right] d\theta dr.
\]
We need the following weighted version of Hardy's inequality:
\begin{lemma} \label{lem:hardy1d}
Let $\psi \in C^2(\RR^+;\RR)$ be such that
\[
\psi(r)\geq 0, \qquad \psi'(r) \leq 0,\qquad
r(\psi'(r)^2 + 2 \psi(r) \psi''(r)) - 2 \psi(r) \psi'(r) \geq 0
\]
for all $r\geq 0$.  Let $f:\RR^+\to\CC$ be such that
 $f(0) = 0$.
Then
\[
\int_0^\infty \psi^2 \frac{|f|^2}{r^2} dr \leq 4 \int_0^\infty \psi^2 |f'|^2 dr
\]
\end{lemma}
{\em Proof:} (inspired by \cite{Sho31}) Let $G$ be the following densely-defined symmetric operator on $L^2(\RR^+)$
\[
G := \frac{1}{i} (\psi\partial_r + \frac{1}{2} \psi')
\]
We have $[G,m] = -i \psi m'$ where $m$ is the operator of multiplication by the function $m(r)$.
We thus have
\begin{eqnarray*}
0 & \leq & \| (G-im) f \|^2 = \langle (G+im)(G-im) f, f\rangle \\
& = & \|Gf\|^2 - \langle (\psi m' - m^2)f,f\rangle
\end{eqnarray*}
Using the definition of $G$,
\begin{eqnarray*}
\| G f\|^2 & = & \int_0^\infty \psi^2|f'|^2 +
\frac{1}{4} (\psi')^2 |f|^2 + \frac{1}{2} \psi\psi' \partial_r|f|^2\ dr \\
& = & \int_0^\infty \psi^2 |f'|^2 +\frac{1}{4} (\psi')^2 |f|^2 -\frac{1}{2} \partial_r (\psi\psi') |f|^2\ dr\\ &&\mbox{} + \frac{1}{2}
\left.\psi(R)\psi'(R)|f(R)|^2\right]_{R=0}^{R=\infty} \\
& \leq & \int_0^\infty \psi^2|f'|^2 -( \frac{1}{4} (\psi')^2 +\frac{1}{2} \psi\psi'')|f|^2 \ dr
\end{eqnarray*}
Thus
\[
0 \leq  \int_0^\infty \psi^2|f'|^2 - (\frac{1}{4} (\psi')^2+\frac{1}{2} \psi\psi'' + \psi m' - m^2)|f|^2\  dr
\]
To establish the lemma we thus need to choose $m$ such that the coefficient of $|f|^2$ in the above
is greater than
$\psi^2/(4r^2)$.  We now check that this is satisfied if we set
$m = -\frac{\psi}{2r}$, provided
\begin{equation}\label{cond:psi}
\frac{1}{4} (\psi')^2+\frac{1}{2} \psi\psi'' -\frac{ \psi \psi'}{2r} \geq 0
\end{equation}
which is equivalent to the condition of the Lemma.
\qed

To apply the above Lemma to $w$, we  set $$f(r) := \left(\int_{|\theta|=1} |w(r,\theta)|^2 d\theta\right)^{1/2}$$ and
\[
\psi(r) := e^{-\sigma r}(1+2\sigma r)^{1/2}.
\]
We check that
\[
\frac{1}{4} (\psi')^2+\frac{1}{2} \psi\psi'' -\frac{ \psi \psi'}{2r} = 3\sigma^4 r^2 e^{-2\sigma r} (1+2r\sigma)^{-1} \geq 0
\]
Moreover,
\begin{eqnarray*}
\int_0^1 \frac{|f(r)|^2}{r^2} dr & = & \int_0^1 \int_{|\theta| = 1} e^{2r\real z}\frac{|u(r\theta)|^2}{r^2} r^{n-1} dr d\theta \\
& \leq & C(z) \|\Omega^{-1} u \|_{L^2(\RR^n)}^2
\end{eqnarray*}
which is finite by Hardy's inequality.  Similarly, $\int_0^1 |f'(r)|^2 dr < \infty$ since $u \in H^1$, and this implies that $f\in C^{1/2}((0,1))$ and thus  $f(0) = 0$.
By the above Lemma then
\begin{equation}\label{whardy}
\int_0^\infty \int_{|\theta|=1} e^{-2r\sigma}(1+2r\sigma)\frac{|w|^2}{4r^2}d\theta dr
\leq \int_0^\infty \int_{|\theta|=1}e^{-2r\sigma}(1+2r\sigma) |\partial_r w|^2 d\theta dr
\end{equation}
Using (\ref{peterpaul}) and  assumption {\bf (A3)}, and taking the limits $R\to \infty,
\epsilon \to 0$, we deduce from (\ref{ibp}) that
\begin{eqnarray*}
\lefteqn{\frac{1}{2}\int_0^\infty\int_{|\theta|=1}e^{-2r\sigma}(1+2r\sigma)
\left\{ (1-a^2)|\partial_r w|^2 + (\delta^2 -
\frac{1}{4}) \frac{|w|^2}{r^2}\right\} d\theta dr}&&\\
\qquad &\leq &\frac{1}{4a^2} \| \Omega f \|_{L^2}^2\hspace{3.5in}
\end{eqnarray*}
and thus using (\ref{whardy}) and optimizing on $a$ obtain
\[
\| e^{-\sigma r} \frac{w}{r}\|_{L^2(\RR^+\times\BS^{n-1})} \leq \frac{1}{2\delta^2} \|
\Omega f \|_{L^2(\RR^n)}
\]
which by the definition of $w$ gives
\[
\| \Omega^{-1} u \|_{L^2(\RR^n)}
\leq \frac{1}{2\delta^2} \| \Omega f \|_{L^2(\RR^n)}
\]
establishing (\ref{est:res}).
\qed
\section{Morawetz and Strichartz estimates}

\subsection{From resolvent to Morawetz}
We recall the  result stating that one can
deduce weighted-$L^2$ spacetime estimates for a Hamiltonian evolution
from a weighted resolvent estimate for the associated elliptic operator
(see Corollary to Theorem XIII.25, \cite[p. 146]{ReeSimIV}.)
\begin{theorem}\label{thm:kato}
(Kato \cite{Kat66}) Let $H$ be a self-adjoint operator on the Hilbert space $X$,  and for $\mu\not\in\RR$
let $\R(\mu) := (H-\mu)^{-1}$ denote the resolvent.
suppose that $A$ is a closed, densely defined operator,
possibly unbounded, from $X$ into a Hilbert space $Y$.  Suppose that
\[
\Gamma := \sup_{\begin{array}{c}\scriptstyle\mu\not\in\RR\\ \scriptstyle\chi\in D(A^*)\\ \scriptstyle\|\chi\|=1\end{array}}
\|A R(\mu) A^* \chi \|_Y < \infty
\]
Then $A$ is $H$-smooth and
\[
\|A\|_H^2 := \sup_{\phi\in X, \|\phi\| = 1} \frac{1}{2\pi} \int_{-\infty}^\infty \| A e^{-itH}\phi\|_Y^2 dt \leq \Gamma^2/\pi^2
\]
\end{theorem}
We use this result to prove (\ref{est:KY}) and (\ref{est:mor}).
\subsubsection{Schr\"odinger's equation}
 Consider
first the case of equation (\ref{eq:lse}).  Set $H=P$, $X = Y = L^2(\RR^n)$ and let
 $A=\Omega^{-1}$, i.e. multiplication by
$\frac{1}{|x|}$.  Thus $A^* = A$, $\R = R$, and if we let $z =  \sqrt{-\mu}$, with the square root branch
 chosen such that $\real z>0$, by Theorem~\ref{thm:ARA*} we then have
\[
\|AR(\mu) A^*\chi\|_Y =
\| \Omega^{-1}(P+z^2)^{-1}\Omega^{-1} \chi
\|_{L^2}  \leq  \frac{1}{2\delta^2} \|\chi\|_{L^2}
\]
Taking the supremum over $\mu$ we see that the hypothesis of Theorem~\ref{thm:kato} is satisfied, and
$\Gamma\leq 1/(2\delta^2)$.
Thus
for $u$ the solution to (\ref{eq:lse}) we have the desired estimate
\[
\| \Omega^{-1} u \|_{L^2_t L^2_x} \leq \frac{1}{\delta^2\sqrt{2\pi}}\|f\|_{L^2}
\]
\subsubsection{The wave equation}
For the wave equation (\ref{eq:lwe}), we instead make the following identifications:
$$X = \Hi^{1/2} \times \Hi^{-1/2},\qquad Y = L^2,\qquad A = (\Omega^{-1}, 0)$$
Recall that $\Hi^s$, defined in Corollary~\ref{corol:equiv}
 are homogeneous Sobolev spaces based on the powers of $P$,
and
thus $$ A^* = (P^{-1/2} \Omega^{-1}, 0). $$ We also let
\[
H = \left(\begin{array}{ll} 0 & -i \\ iP & 0 \end{array}\right)
\]
so that the solution to the wave equation (\ref{eq:lwe}) is
$u = e^{itH}\left(\begin{array}{l} f\\ g\end{array}\right)$.  The resolvent of $H$ is
\[
\R(z) = (H-z)^{-1} = \left( \begin{array}{ll} z(P-z^2)^{-1} & -i(P-z^2)^{-1} \\
 i(P-z^2)^{-1} P & z (P-z^2)^{-1} \end{array}\right)
\]
so that
\[
B:= A \R(z) A^* = \Omega^{-1} z (P-z^2)^{-1} P^{-1/2} \Omega^{-1}
\]
Let
\begin{equation}
D :=  \Omega^{-1} P^{-1} \Omega^{-1}
\end{equation}
\begin{lemma}\label{lem:hardy} The operator $D$ is bounded on $L^2$.
\end{lemma}
{\em Proof:} Let
\[
E := \Omega^{-1} P^{-1/2}
\]
Then $D = EE^*$ and it's thus enough to prove that $E$ is bounded.  This amounts to proving the Hardy inequality
for $P$, i.e.
\[
\| \Omega^{-1} u \|_{L^2} \leq c \| P^{1/2} u \|_{L^2}
\]
which has already been shown (\ref{PHardy}), with $c = 1/\delta$.

We are going to use complex interpolation to prove boundedness of $B$.
This will require the following fact from
operator theory, to be proved in the Appendix.
\begin{theorem}\label{lem:operator_theory}
Suppose $\Lambda$ and $\Omega$ are self-adjoint
operators on a Hilbert space~$X$ and that $\Lambda$ is non-negative
with zero nullspace.
If there is a constant $c$ such that
\begin{equation}
  \label{eq:HH}
    \| [ \Omega , \Lambda ^ 2 ] f \| \le c \| \Lambda f \|
\end{equation}
for all $f \in X$ then $[ \Omega , \Lambda ]$ is a
bounded operator on~$X$.
\end{theorem}

\begin{lemma}\label{lem:pain} The operator $B$ is bounded on $L^2$, uniformly in $z$.
\end{lemma}
{\em Proof:}

We define
\[
  P _ 0 := - \Delta + V _ 0 , \qquad \Lambda := P^{1/2},\qquad \Lambda _ 0 := P _ 0 ^ { 1 / 2 },\qquad R := \Lambda_0\Lambda^{-1}
\]
where
\[
  V _ 0 ( x ) := a |x| ^ { - 2 } , \qquad a \ge 0 .
\]
Let $r = |x|$.  On $\Sigma_l$, the $l$'th spherical harmonic subspace of $L ^ 2 ( \mathbb R ^ n )$, we have
$P_0 = A _ \nu$ where
\[
  A _ \nu := \partial ^ 2 _ r + ( n - 1 ) r ^ { - 1 } \partial _ r
  + \nu ^ 2 r ^ { - 2 }
\]
and
\[
  \nu ^ 2 := ( \lambda + l ) ^ 2 + a , \qquad \lambda := \frac { n - 2 } 2 .
\]
If $a$ is chosen sufficiently large then $\nu ^ 2 > 1$.
We will assume this from now on.  In this case $A _ \nu$ agrees
with the operator of the same name in \cite{PST1} and \cite{BPST1}, the positive branch
having been chosen for the square root.
For $n > 4$ we can choose $a = 0$, but this does not result in any
real simplification of the argument below.

Let us define the following operators:
\[
  E_ 0 := \Omega ^ { - 1 } \Lambda _ 0 ^ { - 1 } ,\quad M_0 := [ \Lambda _ 0 ^ 2 , \Omega ] \Lambda _ 0 ^ { - 1 },
  \quad
  C _ 1 := \Omega \Lambda _ 0 ^ { - 1 } \Omega ^ { - 1 } \Lambda _ 0 ,\quad
  C _ 2 := \Omega \Lambda _ 0 \Omega ^ { - 1 } \Lambda _ 0 ^ { - 1 }.
\]
Note that $E_0$ is the same as the operator $E = \Omega^{-1}\Lambda^{-1}$, but with the potential $V$ replaced by $V_0$, hence the $L^2$-boundedness of $E_0$ follows from Hardy's inequality by the same argument as in the proof of Lemma~\ref{lem:hardy}.  The boundedness of $M_0$ is reduced to that of $E_0$ by a direct computation.  Meanwhile,
$$ C_2 = [\Omega,\Lambda_0]E_0 + I.$$
Thus the boundedness of $C_2$ follows from that of $M_0$, using Theorem~\ref{lem:operator_theory}. Finally we have the following lemma, the proof of which will be given in the Appendix.
\begin{lemma}\label{lem:mellin}
The operator $C_1$ is bounded on $L^2$.
\end{lemma}

To proceed with complex interpolation, we define the following family of operators: For $s\in \mathbb C$, $0\leq \Re s \leq 1$, let
\[
  T _ s := z ^ { 2 s } e ^ { s ^ 2 - 1 / 4 } \Omega ^ { - 1 }
    ( P - z ^ 2 ) ^ { - 1 } \Lambda _ 0 ^ { - 1 } \Omega ^ { - 1 }
    \Lambda _ 0 ^ { 1 - 2 s } .
\]
Up to a constant,
\[
  T _ 0 =  \Omega ^ { - 1 } ( P - z ^ 2 ) ^ { - 1 } \Lambda _ 0 ^ { - 1 }
    \Omega ^ { - 1 } \Lambda_0
  = \Omega ^ { - 1 } ( P - z ^ 2 ) ^ { - 1 } \Omega ^ { - 1 } C _ 1 .
\]
Meanwhile
$$
  T _ 1 =  \Omega ^ { - 1 } z ^ 2 ( P - z ^ 2 ) ^ { - 1}
    \Lambda _ 0 ^ { - 1 } C _ 0
    = \Omega ^ { - 1 } ( P ( P - z ^ 2 ) ^ { - 1 } - I )
    \Lambda _ 0 ^ { - 1} C _ 0
    = \Omega ^ { - 1 } ( P - z ^ 2 )^{-1} \Omega ^ { - 1 } F
      - C _ 0 ^ 2,
$$
where
\[
  F := \Omega P \Lambda _ 0 ^ { - 1 } C _ 0
    = C _ 2 + ( \Omega ^ 2 V - a ) C _ 0 ^ 2.
\]
By {\bf (A1)},
$\Omega ^ 2 V$ is bounded in $L ^ \infty$.
We have proved the uniform $L ^ 2$ boundedness of
$\Omega ^ { - 1 } ( P - z ^ 2 )^{-1} \Omega ^ { - 1 }$ for $z^2\not\in \mathbb{R}^+$ in Theorem~\ref{thm:ARA*}.
So $T _ 0$ and $T _ 1$ are bounded. The contribution of the imaginary
part of $s$ to the power of $\Lambda _ 0$ only puts a unitary
operator at the tail end, which doesn't affect boundedness,
and its contribution to the $z$ power is taken care of by the exponential term.  Therefore $T_{0+it}$ and $T_{1+it}$ are bounded uniformly in $t$ and in $z$.  By complex interpolation $T _ { 1 / 2 }$ is bounded uniformly in $z$:
\[
  T _ { 1 / 2 } = z \Omega ^ { - 1 } ( P - z ^ 2 ) ^ { - 1 }
    \Lambda _ 0 ^ { - 1 } \Omega ^ { - 1 }.
\]
We have
\[
  B = T _ { 1 / 2 } \Omega \Lambda _ 0 \Lambda ^ { - 1 } \Omega ^ { - 1 }
    = :T _ { 1 / 2 } G
\]
and
\[
  G = \Lambda _ 0 \Omega \Lambda ^ { - 1 } \Omega ^ { - 1 }
    + [ \Omega , \Lambda _ 0 ] \Lambda ^ { - 1 } \Omega ^ { - 1 } =: J + [\Omega,\Lambda_0] E^*
\]
Finally,
\[
  J = \Lambda _ 0 \Lambda ^ { - 1 } + \Lambda _ 0 [ \Omega, \Lambda ^ { - 1 } ]
    \Omega ^ { - 1 }
  = R + R [ \Lambda , \Omega ] E^*.
\]
 Moreover,
\[
 M := [ \Lambda ^ 2 , \Omega ] \Lambda^{-1}= [ \Lambda _ 0 ^ 2 , \Omega ] \Lambda^{-1} = M_0 R,
\]
Note that $R$ is bounded by the equivalence
of Sobolev norms (Corollary~\ref{corol:equiv}). Therefore $M$ is bounded too.
The boundedness of $[ \Lambda , \Omega ]$ and $[ \Lambda _ 0 , \Omega _ 0 ]$
now follows from Theorem~\ref{lem:operator_theory}.  These in turn imply the boundedness of $J$, $G$, and hence that of $B$.  This concludes the proof of Lemma~\ref{lem:pain}.

We can therefore
again apply Theorem~\ref{thm:kato} to deduce the Morawetz estimate
(\ref{est:mor}), except that the norms on the right will be $P$-based norms. By the equivalence
result of Corollary~\ref{corol:equiv} however, we can replace those with standard Sobolev norms.

\subsection{From Morawetz to Strichartz}
\subsubsection{Schr\"odinger's equation}
We consider the
potential term as a source term,
\begin{equation}
  \label{eq:ff1}
  i\partial_t u+\Delta u=V(x) u,\qquad u(0) = f
\end{equation}
and integrate using $S_0(t) =e^{-it\Delta}$, the free evolution, to get
\begin{equation}
  \label{eq:duhamel}
  u(t) =S_0(t)f+\int_0^t S_0(t-s) V u(s) ds
\end{equation}
 The first term can be ignored  since it satisfies the estimate we want to prove,
and we can focus on the Duhamel term. Given that $n\geq 3$,  one has
Strichartz estimates up to the end-point for the free
evolution, \cite{KeeTao98}, i.e. for the
pair $(p,q)=(2,\frac{2n}{n-2})$. We recall that these Strichartz
estimates hold in a slightly relaxed setting,
\begin{equation}
  \label{eq:yeah}
\|  \int_0^t S_0(t-s) F(x,s)ds \|_{L^2_t(L_x^{\frac{2n}{n-2},2})} \leq C \| F \|_{L^2_t(L_x^{\frac{2n}{n+2},2})} ,
\end{equation}
where $L^{\alpha,\beta}$ are Lorentz spaces. Hence to prove our
estimate, all we need to check is $F = Vu \in
L^2_t(L_x^{\frac{2n}{n+2},2})$. However, from (\ref{est:KY}) we
have $\Omega^{-1} u\in L^2_t L^2_x$, while assumption {\bf (A1)} implies
$\Omega V \in L^{n,\infty}$. Thus, using O'Neil's inequality (H\"older
inequality for Lorentz spaces \cite{One63}) we have
\begin{eqnarray*}
\|  \int_0^t S_0(t-s) V u(s) ds \|_{L^2_t(L_x^{\frac{2n}{n-2},2})} &
\leq & C \| V u \|_{L^2_t(L_x^{\frac{2n}{n+2},2})}\\
& \leq & C \| \Omega V \|_{L^{n,\infty}} \| \Omega^{-1} u\|_{L^2_tL^2_x}
\\
& \leq & C \|f \|_{L^2}
\end{eqnarray*}
which proves (\ref{eq:strichartz}) at the end-point
$(p,q)=(2,\frac{2n}{n-2})$.  Interpolating between this and the
conservation of the $L^2$ norm for (\ref{eq:lse}), which corresponds to
$(p,q) = (\infty,2)$ in (\ref{eq:strichartz}), one obtains the full
range of Strichartz estimates.

\subsubsection{Wave equation}

We write the solution to (\ref{eq:lwe}) as the sum
of the  solution to the free wave equation plus a Duhamel term
\begin{equation}
  \label{eq:ff10}
  u(t)=\dot{W}(t)f+W(t)g-\int_0^t W(t-s) V(x) u(s)ds,
\end{equation}
where $W(t)=\frac{\sin(t\sqrt{-\Delta})}{\sqrt{-\Delta}}$, and
$\dot{W}=\partial_t W$. We again ignore the first two terms in the above and focus on the Duhamel
term.  Since $W(t-s) = -\dot{W}(t) {W}(s) + {W}(t)\dot{W}(s)$, this splits into
two terms.  We will deal with the
first one, the treatment of the second term being similar. We are going to use
the following lemma,
\begin{lemma}[\cite{ChrKis01}]
\label{lem:CK}
 Let $X,Y$ be two Banach spaces and let $T$ be a bounded linear operator from $L^\beta(\RR^+;X)$ to
 $L^\gamma(\RR^+;Y)$, $Tf(t)=\int_0^\infty K(t,s)f(s)ds$. Then the
 operator $\tilde{T} f(t)=\int_0^t K(t,s)f(s)ds$ is bounded from
 $L^\beta(\RR^+;X)$ to $L^\gamma(\RR^+;Y)$ when $\beta<\gamma$, and $\|\tilde{T}\| \leq c_{\beta,\gamma} \|T\|$ with
 $c_{\beta,\gamma} = (1-2^{1/\gamma-1/\beta})^{-1}$.
\end{lemma}
We set
$$
T h(t) :=\dot{W}(t)\int W(s) \Omega^{-1} h(s) ds.
$$
Using the following Strichartz estimate for the free wave equation
\[
\|\dot W(t) F \|_{L^p_x\dot{H}^\sigma_q} \leq C \|F\|_{\dot{H}^{1/2}(\RR^n)},
\]
combined with the {\em dual} to the Morawetz
estimate (\ref{est:mor}) for the free wave equation, namely
\[
\| \int W(s) G(s) ds \|_{\dot{H}^{1/2}(\RR^n)} \leq C \|\Omega G\|_{L^2(\RR^{n+1})}
\]
we obtain
\begin{eqnarray*}
\|Th\|_{L^p\dot{H}^\sigma_q} & \leq & C \|\int W(s) \Omega^{-1} h(s) ds\|_{\dot{H}^{1/2}(\RR^n)}\\
& \leq &C \|h\|_{L^2(\RR^{n+1})}
\end{eqnarray*}
with $p$, $q$ and $\sigma$ as in the statement of the Theorem.  By Lemma~\ref{lem:CK}, the corresponding operator
$\tilde{T}$ satisfies the same estimate as $T$ (with a different constant).  On the other hand, the solution to
(\ref{eq:lwe}) is
\[
u(t) = \dot{W}(t)f+W(t)g +  \tilde{T}(\Omega V u)
\]
By assumption {\bf (A1)} and (\ref{est:mor}) we conclude
\begin{eqnarray*}
\|\Omega V u \|_{L^2(\RR^{n+1})} & \leq & \max\{\gamma_+^2,\gamma_-^2\} \| \Omega^{-1} u \|_{L^2(\RR^{n+1})} \\
& \leq & C ( \| f\|_{\dot{H}^{1/2}} + \|g \|_{\dot{H}^{-1/2}})
\end{eqnarray*}
which establishes (\ref{eq:strichartzwave}).

\section{The point-dipole potential}
An example of a physical potential  satisfying our assumptions {\bf (A1-A3)}
is that of an electrical point-dipole.  The Schr\"odinger equation with this
potential arises
for example in the study of electron capture by polar molecules \cite{Lev67}. Let $\psi$
be the wave function of an electron in the electric field of a dipole that is supposed to be point-like
and fixed at the origin.  The equation then reads
\begin{equation}\label{eq:elecdip}
i \partial_t \psi = - \Delta \psi + \frac{\mathbf{p}\cdot x}{|x|^3} \psi
\end{equation}
where $\mathbf{p} := \frac{2me\mathbf{D}}{\hbar^2}$ is dimensionless.  Here $m,e$
are the mass and charge of the electron and $\mathbf{D}$ is the electric dipole
moment of the molecule.  Choosing coordinates such that $\mathbf{p} =
(0,0,p)$,
The potential $V(x) = px_3/|x|^3$ is homogeneous of degree -2, so that assumptions {\bf (A2)} and
{\bf (A3)} coincide, and for the weighted-$L^2$ (\ref{est:KY}) and Strichartz (\ref{eq:strichartz}) estimates
to hold for $\psi$, all we need is that the lowest eigenvalue of the
operator $-\sphlap + p x_3$ on $\BS^2$ be larger than $-1/4$.
 This is clearly the case if $p<1/4$, and if we
let $p_0$ denote the largest value of $p$ for which
this continues to
hold, it is known that $p_0\approx 1.28$ (see \cite{Lev67} for the calculation of this
``critical value" of the dipole moment).

\section{Appendix}
\subsection{Proof of Theorem~\ref{lem:operator_theory}}
From this point on $c$ will denote a constant which depends only
on the dimension, but whose value may differ from equation
to equation.

We define rescaled versions of $\Omega$ and $\Lambda$, $$\Omega _ s := s
^ { - 1 / 2 } \Omega,\qquad \Lambda _ s := s ^ { 1 / 2 } \Lambda.$$
We have
\begin{lemma}\label{lem:On_Q} For $\alpha \geq 0$, let $ Q _ \alpha ( \sigma ) := \Lambda _
  s ^ { \alpha } \exp ( - \Lambda _ s ^ 2 )$, then
\begin{enumerate}
\item the operator $Q_\alpha(\sigma)$ is bounded on $L^2$, uniformly in
  $s$, for all $\alpha\geq 0$.
\item For $\alpha, \gamma > 0$,
  \begin{equation}
    \label{eq:RH}
     \Lambda ^ { \alpha - 2 \gamma }= \Gamma(\gamma)^{-1} \int _ 0 ^ \infty s ^ { \gamma - \alpha / 2 } Q _ \alpha ( s )
    \frac { d s } { s }.
  \end{equation}
\item For $\alpha > 0$,
  \begin{equation}
    \label{eq:L2}
      \int _ 0 ^ \infty \| Q _ \alpha ( s ) g \| ^ 2 \frac { d s } s
  = \displaystyle 2 ^ { - \alpha } \Gamma ( \alpha ) \| g \| ^ 2 .
  \end{equation}
\end{enumerate}
\end{lemma}

By spectral theory we may assume without loss of generality
that $X = L ^ 2 ( \mathcal X )$ where $\mathcal X$ is a
measure space and that $\Lambda _ s$ is multiplication by
a non-negative real measurable function on $\mathcal X$,
\[
   ( \Lambda
g ) ( \xi  ) = m ( \xi  ) g ( \xi  )\,.
\]
All of the results above then follow immediately.  For example,
\[\begin{array}{r@{\:=\:}l}
  \displaystyle \left ( \int _ 0 ^ \infty s ^ { \gamma - \alpha / 2 }
    Q _ \alpha ( s ) \frac { d s } { s } g \right ) ( \xi  )
  & \displaystyle \int _ 0 ^ \infty s ^ { \gamma } m^\alpha ( \xi  )
    \exp ( - s m^2 ( \xi  ) ) \frac { d s } { s } g ( \xi  )
  \cr & \displaystyle \Gamma ( \gamma )
    m ^ { \alpha - 2 \gamma } ( \xi ) g ( \xi  )
  \cr & \displaystyle \Gamma ( \gamma )
    \left ( \Lambda ^ { \alpha - 2 \gamma } g \right ) ( \xi  ),
\end{array}\]
so that
\[
  \int _ 0 ^ \infty s ^ { \gamma - \alpha / 2 } Q _ \alpha ( s )
    \frac { d s } { s }
  = \Gamma ( \gamma ) \Lambda ^ { \alpha - 2 \gamma } .
\]
Similarly,
\[\begin{array}{r@{\:=\:}l}
  \displaystyle
  \int _ 0 ^ \infty \| Q _ \alpha ( s ) g \| ^ 2 \frac { d s } s
 & \displaystyle
  \int _ 0 ^ \infty \int _ { \xi  \in \mathcal X }
  s ^ { \alpha } m ^ { 2 \alpha } ( \xi ) \exp ( - 2 s m ^2 ( \xi  ) ) | g ( \xi ) | ^ 2
  \, d \xi  \, \frac { d s } s
  \cr & \displaystyle
  \int _ { \xi  \in \mathcal X } \int _ 0 ^ \infty
  s ^ \alpha m ^ { 2 \alpha } ( \xi ) \exp ( - 2 s m ^2 ( \xi  ) )
  \, \frac { d s } s  | g ( \xi  ) | ^ 2 \, d \xi
  \cr & \displaystyle
  \int _ { \xi  \in \mathcal X } 2 ^ { - \alpha } \Gamma ( \alpha ) | g ( \xi  ) | ^ 2 \, d \xi
  \cr & \displaystyle 2 ^ { - \alpha } \Gamma ( \alpha ) \| g \| ^ 2 .
\end{array}\]

For our present purposes, the main use of the operators $Q _ \alpha (
s )$ is the following boundedness criterion, which may be seen as a
simple form of the Cotlar-Stein lemma.
\begin{lemma}\label{lem:localisation}
Let $\epsilon>0$, and suppose that $\mathcal{T}$ is an operator such that
\begin{equation}
  \label{eq:CS}
  \| Q _ \alpha ( s ) \mathcal{T} Q _ \alpha ( t ) \| \le c \exp ( -
  \epsilon|\log s -\log t| ) ,
\end{equation}
then
\[
  \| \mathcal{T} \| \le 2 ^ { 3 / 2 } \epsilon ^ { - 3 / 2 } \Gamma ( \alpha ) ^ { - 1 } c .
\]
\end{lemma}
 By Lemma~\ref{lem:On_Q} we may write
\[
  Q _ \alpha ( s ) \mathcal{T} = \Gamma ( \alpha ) ^ { - 1 } \int _ 0 ^ \infty Q _ \alpha ( s ) \mathcal{T} Q _ { 2 \alpha } ( 2 t ) \frac { d t } t ,
\]
and, by the triangle inequality,
\[
  \| Q _ \alpha ( s ) \mathcal{T} f \| \le \Gamma ( \alpha ) ^ { - 1 } \int _ 0 ^ \infty \| Q _ \alpha ( s ) \mathcal{T} Q _ { \alpha } ( t ) \| \| Q _ { \alpha } ( t ) f \| \frac { d t } t .
\]
Let $$d(s,t)=e^{ |\log s-\log t| }.$$ Using \eqref{eq:CS} and then Cauchy-Schwarz,
\[\begin{array}{r@{\:}c@{\:}l}
  \displaystyle \| Q _ \alpha ( s ) \mathcal{T} f \|
 & \le & c \Gamma ( \alpha ) ^ { - 1 } \int _ 0 ^ \infty d ( s , t ) ^ { - \epsilon } \| Q _ { \alpha } ( t ) f \| \frac { d t } t \cr
  & \le & \displaystyle
  c \Gamma ( \alpha ) ^ { - 1 }
    \left [ \int _ 0 ^ \infty d ( s , t ) ^ { - \epsilon } \frac { d t } t \right ] ^ { 1 / 2 }
    \left [ \int _ 0 ^ \infty d ( s , t ) ^ { - \epsilon } \| Q _ { \alpha } ( t ) f \| ^ 2 \frac { d t } t \right ] ^ { 1 / 2 }
  \cr & = & \displaystyle 2 \epsilon ^ { - 1 } c \Gamma ( \alpha ) ^ { - 1 }
    \left [ \int _ 0 ^ \infty d ( s , t ) ^ { - \epsilon } \| Q _ { \alpha } ( t ) f \| ^ 2 \frac { d t } t \right ] ^ { 1 / 2 } .
\end{array}\]
Then, squaring this last inequality and integrating over $s$,
\[
  \int _ 0 ^ \infty \| Q _ \alpha ( s ) \mathcal{T} f \| ^ 2 \frac { d s } s
  \le 4 \epsilon ^ { - 2 } c ^ 2 \Gamma ( \alpha ) ^ { - 2 } \int _ 0 ^ \infty \int _ 0 ^ \infty d ( s , t ) ^ { - \epsilon } \| Q _ { \alpha } ( t ) f \| ^ 2 \frac { d t } t \frac { d s } s
\]
or, switching the order of integration on the right and using \eqref{eq:L2},
\[
  \| \mathcal{T} f \| ^ 2 \le 8 \epsilon ^ { - 3 } c ^ 2 \Gamma ( \alpha ) ^ { - 2 } \| f \| ^ 2
\]
which is the claim.
\begin{lemma}\label{lem:On_K}
Define $ K_\alpha ( \sigma ) := [ \Omega _ s , Q _ \alpha(
\sigma ) ]$ and $ L ( s ):= [ \Lambda _ s ^ 2 , \Omega _ s ] Q _ 0 ( s
)$.
\begin{enumerate}
\item Both operators $K _ 0 ( s )$ and $L(s)$ are bounded on $L ^ 2$, uniformly in $s$.
\item The operator $Q_2 (r) K_2(s)$ is bounded on $L ^ 2$, and
  \begin{equation}
    \label{eq:QK}
    \| Q _ 2 ( r ) K _ 2 ( s ) \| \le c d ( r , s ) ^ { - 1/2 }\,.
  \end{equation}
\end{enumerate}
\end{lemma}
We start with $K_0$~:
$$
  \partial _ s ( s ^ { 1 / 2 } K _ 0 ( s ) f)
     =  [ \Omega , \partial_s  Q _ 0 ( s ) ]f = [ \Omega , - \Lambda ^ 2 Q _ 0 ( s ) ] f
     =  - \Lambda ^ 2 ( s ^ { 1 / 2 } K _ 0 ( s ) )
    f + [ \Lambda ^ 2 , \Omega ] Q _ 0 ( s )f\, .
$$
Thus, taking the scalar product with $ s ^ { 1 / 2 } K _ 0 ( s ) f$,
\begin{equation*}
 \partial _ s \| s ^ { 1 / 2 } K _ 0 ( s ) f \| _ 2
   \leq  2 \| [ \Lambda ^ 2 , \Omega ] Q _ 0 ( s ) f \| _ 2
   \leq  c s ^ { - 1 / 2 } \| Q _ 1 ( s ) f \| _ 2
  \leq c s ^ { - 1 / 2} \| f \| _ 2 ,
\end{equation*}
where we have used that $\Lambda$ satisfies \eqref{eq:HH}.  Integrating both sides of
the above on $[ 0 , s ]$ we obtain the boundedness of $K _ 0$.

Next, using again \eqref{eq:HH},
\[\begin{array}{r@{\:}c@{\:}l}
  \| L ( s ) f \|
  & = & \| [ \Lambda _ s ^ 2 , \Omega _ s ] Q _ 0 ( s ) f \|
   =  \| s ^ { - 1 / 2 } [ \Lambda ^ 2 , \Omega ] Q _ 0 ( s ) f \|
  \cr & \le & c \| s ^ { - 1 / 2 } \Lambda  Q _ 0 ( s ) f \|
   =  c \| Q _ 1 ( s ) f \|\le  c \| f \| .
\end{array}\]
Finally, we establish the bound on $Q _ 2 ( r ) K _ 2 ( s )$.
We can write it in either of the two forms
\begin{eqnarray*}
  Q_2 (r) K_2(s) & = & \sqrt \frac s r Q_2(r)
    [\Omega_r,\Lambda_r^2] Q_0(s)+Q_2(r)\Lambda_s^2 \Omega_s
    Q_0(s)-Q_2(r) Q_2(s) \Omega_s\\
 & = & \sqrt \frac s r Q_2(r)
    [\Omega_r,\Lambda_r^2] Q_0(s)+\frac s r Q_4(r)K_0(s)\\
 & = & 2^{3/2} \sqrt \frac s r Q_2(r/2)(-L^*(r/2)) Q_0(s)+\frac s r Q_4(r)K_0(s).
\end{eqnarray*}
or
\begin{eqnarray*}
    Q _ 2 ( r ) K _ 2 ( s )
  & = & \frac r s Q _ 0 ( r )
  \left (
    [ \Lambda _ s ^ 2 , \Omega ] Q_2(s)
    + [ \Omega _ s , Q _ 4 ( s ) ]
  \right ) \\
 & = & 2 ^ { 3 / 2 } \frac r s Q _ 0 ( r ) \left (
    L ( s/2 ) Q _ 2 ( s/2 )
    + K _ 2 ( s/2 ) Q _ 2 ( s/2 )
    + Q _ 2 ( s/2 ) K _ 2 ( s/2 )
  \right ) \\
 & = & 2 ^ { 3 / 2 } \frac r s Q _ 0 ( r ) \left (
    L ( s/2 ) Q _ 2 ( s/2 )
    - L ( s/2 ) ^ \star Q _ 2 ( s/2 )
    + K _ 0 ( s/2 ) Q _ 4 ( s/2 )
   \right . \\
 & & \qquad \left . {}
    - Q _ 2 ( s/2 ) L ( s/2 )
    + Q _ 4 ( s/2 ) K _ 0 ( s/2 )
  \right )
\end{eqnarray*}
From the uniform boundedness of $Q _ \alpha$, $K _ 0$ and $L$
we see that
\[
  \| Q _ 2 ( r ) K _ 2 ( s ) \| \le c \sqrt \frac s  r \,\,\text{   and    }\,\,
  \| Q _ 2 ( r ) K _ 2 ( s ) \| \le c \sqrt \frac r  s ,
\]
which yields \eqref{eq:QK}.

Theorem~\ref{lem:operator_theory} will follow from the following identity:
$$
  [\Lambda,\Omega] =  \int_0^\infty Q_1(s) \Lambda \Omega\,\frac{ds}
  s -\int_0^\infty \Omega \Lambda  Q_1(s) \,\frac {ds}s
  =   \int_0^\infty Q_2(s) \Omega_s -\Omega_s  Q_2(s) \,\frac {ds}s
 =   \int_0^\infty K_2(s) \,\frac {ds}s.
$$
The forthcoming Lemma ~\ref{lem:E_st} then provides a bound on
$Q_2(r)[\Lambda,\Omega]Q_2(s)$ which we use to apply
Lemma~\ref{lem:localisation} and obtain boundedness of $[\Lambda,\Omega]$.
\begin{lemma}\label{lem:E_st}
Define
\[
  E ( r , t ) := \int _ 0 ^ \infty Q _ 2 ( r ) K _ 2 ( s ) Q _ 2 ( t ) \frac { d s } s .
\]
Then
\[
  \| E ( r , t ) \| \le c d ( r , t ) ^ { - 1 / 4 }.
\]

\end{lemma}
Since $E ( r , t ) ^ * = - E ( t , r )$ we may assume without
loss of generality that $r \le t$.
We write
\[
  E ( r , t ) = E _ < ( r , t ) + E _ > ( r , t )
\]
where
\[
   E _ < ( r , t ) := \int _ 0 ^ { \sqrt { r t } } Q _ 2 ( r ) K _ 2 ( s ) Q _ 2 ( t ) \frac { d s } s,\qquad
 E _ > ( r , t ) := \int _ { \sqrt { r t } } ^ \infty Q _ 2 ( r ) K _ 2 ( s ) Q _ 2 ( t ) \frac { d s } s.
\]
Remark that $K_2(s) Q_2(t) =- (Q_2(t)K_2(s))^*$. Then, we use
\eqref{eq:QK} and for $s \le \sqrt { r t }$,
\[
  \| Q _ 2 ( r ) K _ 2 ( s ) Q _ 2 ( t ) \| \le \| Q _ 2 ( r ) \| \| K _ 2 ( s ) Q _ 2 ( t ) \| \le c d ( s , t ) ^ { - 1/2 } = c \sqrt{s / t} .
\]
while for $s \ge \sqrt { r t }$,
\[
  \| Q _ 2 ( r ) K _ 2 ( s ) Q _ 2 ( t ) \| \le \| Q _ 2 ( r ) K _ 2 ( s ) \| \| Q _ 2 ( t ) \| \le c d ( r , s ) ^ { - 1/2 } = c \sqrt{r / s} .
\]
Integrating,
\[
  \| E _ < ( r , t ) \| \le c ( r / t ) ^ { 1/4 } = c d ( r , t ) ^ { - 1 / 4 }
\]
and
\[
  \| E _ > ( r , t ) \| \le c ( r / t ) ^ { 1 /4 } = c d ( r , t ) ^ { - 1 / 4 } ,
\]
which completes the proof.
\subsection{Proof of Lemma~\ref{lem:mellin}}
Recall the following  from \cite{PST1}:
The Mellin transform, its inverse
\begin{displaymath}
(\mathcal{M}\phi)(z) = \int^\infty_0 r^{z-n} \phi(r) r^{n-1}\,dr,
\qquad
(\mathcal{M}^{-1}f)(r) = -\frac{1}{2\pi i}\int_C r^{-z}f(z)\,dz,
\end{displaymath}
and its action on multiplication by powers of $r$,
\begin{displaymath}
(\mathcal{M}\Omega^\sigma\phi)(z) = (\mathcal{M}\phi)(z+\sigma).
\end{displaymath}
We also recall the Hankel transform,
\begin{displaymath}
(\mathcal{H}_\nu\phi)(r) = \int^\infty_0(rs)^{-\lambda}J_\nu(rs)\phi(s)s^{n-1}\,ds,
\end{displaymath}
where $J_\nu$ is the Bessel function of the first kind of order $\nu$,
its composition with the Mellin transform,
\begin{displaymath}
(\mathcal{M}\mathcal{H}_\nu\phi)(z) = 2^{z-\lambda-1}
\frac{\Gamma(\frac{z-\lambda+\nu}{2})}{\Gamma(1-\frac{z-\lambda-\nu}{2})}
(\mathcal{M}\phi)(2\lambda+2-z)
\end{displaymath}
and the representation of $\Lambda_0^\sigma$ on $\Sigma_l$ via the Hankel transform
\begin{displaymath}
(A_\nu)^{\sigma/2} = \mathcal{H}_\nu\Omega^\sigma\mathcal{H}_\nu.
\end{displaymath}
Now Plancherel's formula can be written in the form
\begin{displaymath}
\left<\varphi,\psi\right> = \frac{\mathrm{Vol}(S^{n-1})}{2\pi}
\int_{-\infty}^\infty(\mathcal{M}\varphi)(\lambda+1+iy)\overline{(\mathcal{M}\psi)(\lambda+1+iy)}\,dy,
\end{displaymath}
from which
\begin{displaymath}
\|\varphi\|^2 = \frac{\mathrm{Vol}(S^{n-1})}{2\pi}
\int_{-\infty}^\infty|(\mathcal{M}\varphi)(\lambda+1+iy)|^2\,dy.
\end{displaymath}
Hence, if $\mathcal O$ is an operator whose action is given in terms of the Mellin transform by
\begin{displaymath}
(\mathcal{M}\mathcal{O}\varphi)(z)=O(z)(\mathcal{M}\varphi)(z)
\end{displaymath}
then its $L^2$ to $L^2$ norm is
\begin{displaymath}
\|\mathcal{O}\| = \sup|O(\lambda+1+iy)|.
\end{displaymath}
Applying this with $\mathcal{O} = \left. C_1\right|_{\Sigma_l}$ we obtain
\begin{displaymath}
\mathcal{O} = \Omega A_\nu^{-1/2}\Omega^{-1}A_\nu^{1/2} = \Omega\mathcal{H}_\nu
\Omega^{-1}\mathcal{H}_\nu\Omega^{-1}\mathcal{H}_\nu\Omega\mathcal{H}_\nu ,
\end{displaymath}
whose action on the Mellin transform side is multiplication by
\begin{displaymath}
O(z) =
\frac{\Gamma(\frac{\nu+z-\lambda+1}2)^2}{\Gamma(\frac{\nu-z+\lambda+1}2)^2}
\frac{\Gamma(\frac{\nu-z+\lambda}2)}{\Gamma(\frac{\nu+z-\lambda}2)}
\frac{\Gamma(\frac{\nu-z+\lambda+2}2)}{\Gamma(\frac{\nu+z-\lambda+2}2)}.
\end{displaymath}
We thus get
\begin{displaymath}
\|\mathcal{O}\| = \frac {\nu^2}{\nu^2-1}
\end{displaymath}
for $\nu > 1$, because
\begin{displaymath}
|O(\lambda+1+iy)| = \left|
\frac{\Gamma(\frac{\nu+iy+2}2)^2}{\Gamma(\frac{\nu-iy}2)^2}
\frac{\Gamma(\frac{\nu-iy-1}2)}{\Gamma(\frac{\nu+iy+1}2)}
\frac{\Gamma(\frac{\nu-iy+1}2)}{\Gamma(\frac{\nu+iy+3}2)}
\right| = \left| \frac { ( \nu + i y ) ^ 2 } { ( \nu + i y ) ^ 2 - 1 } \right|.
\end{displaymath}
For the $P _ 0$ on the $l$'th spherical harmonic subspace
we have
\[
  \nu = \sqrt { ( \lambda + l ) ^ 2 + a } \ge \sqrt { \lambda ^ 2 + a } > 1,
\]
for $a>1-\lambda^2$.  This show the boundedness of $\left. C _ 1\right|_{\Sigma_l}$.  The $L^2$ boundedness of $C_1$ then follows from the orthogonality of the subspaces $\Sigma_l$ with respect to the $L^2$ innerproduct.

\bigskip

{\bf Acknowledgements} We are grateful to Piero D'Ancona for pointing out an error in the original proof we had of Lemma \ref{lem:pain}.
\bibliographystyle{plain}
\def\cprime{$'$}

\end{document}